\documentclass[fleqn]{mat01}
\usepackage{times,mathtimy,amssymb,latexsym}
\begin{document}

\setcounter{page}{147}
\firstpage{147}

\newtheorem{theore}{Theorem}
\renewcommand\thetheore{\arabic{section}.\arabic{theore}}
\newtheorem{theor}{\bf Theorem}
\newtheorem{lem}[defin]{\it Lemma}
\newtheorem{propo}[theore]{\rm PROPOSITION}
\newtheorem{coro}[theore]{\rm COROLLARY}
\newtheorem{definit}[theore]{\rm DEFINITION}
\newtheorem{probl}[theore]{\it Problem}
\newtheorem{exampl}[theore]{\it Example}
\newtheorem{pot}{\it Proof of Theorem}
\newtheorem{poc}[theore]{\it Proof of Corollary}
\newtheorem{remar}[theore]{\it Remark}
\newtheorem{conjec}{\it Conjecture}

\def\deff{\trivlist \item[\hskip \labelsep{DEFINITION.}]}
\def\probl{\trivlist \item[\hskip \labelsep{\bf Problem.}]}
\def\claim{\trivlist \item[\hskip \labelsep{\it Claim.}]}

{\def\nequiv{\equiv\hskip -.7pc\lower 1pt\hbox{/}\hskip .7pc}

\title{A variant of Davenport's constant}

\markboth{R Thangadurai}{A variant of Davenport's constant}

\author{R THANGADURAI}

\address{Harish-Chandra Research Institute, Chhatnag Road,
Jhusi,  Allahabad 211~019, India\\
\noindent E-mail: thanga@hri.res.in}

\volume{117}

\mon{May}

\parts{2}

\pubyear{2007}

\Date{MS received 25 July 2006; revised 27 September 2006}

\begin{abstract}
Let $p$ be a prime number. Let $G$ be a finite abelian
$p$-group of exponent $n$ (written additively) and $A$ be a
non-empty  subset   of $]n[:=  \{1,2,\dots, n\}$ such that
elements of $A$ are incongruent modulo $p$ and non-zero modulo $p$.
Let $k \geq D(G)/|A|$ be any integer where $D(G)$ denotes the
well-known Davenport's constant. In this article, we prove that
for any sequence $g_1, g_2, \dots, g_k$ (not necessarily
distinct) in  $G$, one can always extract a subsequence $g_{i_1},
g_{i_2}, \dots, g_{i_\ell}$ with  $1\leq  \ell \leq k$ such that
\begin{equation*}
\sum_{j=1}^\ell a_{j}g_{i_j}  = 0 \mbox{ in } G,
\end{equation*}
where $a_j \in A$ for all $j$. We provide examples where this
bound cannot be improved. Furthermore, for the cyclic groups,
we prove some sharp results in this direction. In the last
section, we explore  the relation between this problem and a
similar problem with prescribed length. The proof of Theorem~1 uses
group-algebra techniques, while for the other theorems, we use
elementary number theory techniques.
\end{abstract}

\keyword{Davenport's constant; zero-sum problems; abelian groups.}

\maketitle

\section{Introduction}

Let $G$ be a finite abelian group additively written. Let $n$ be
the exponent of $G$. Let $\emptyset \ne$ \hbox{$A\subset\ ]n[$}
where $]n[ := \{1,2,\dots,n\}$. The set of all integers is denoted
by ${\Bbb Z}$, while the set of all positive integers is denoted
by ${\Bbb N}$.  Also, let $p$ be a prime number. The finite field
with $p$ elements is denoted by ${\Bbb F}_p$ or ${\Bbb Z}_p$.
Also, direct sum of $d$ copies of ${\Bbb Z}_p$ is denoted by
${\Bbb Z}_p^d$. For any integer $x\ge 1$, we denote $]x[$ for
$\{1,2,\dots,x\}$.

Since $G$ is an abelian group, $G$ is a ${\Bbb Z}$-module. As a
${\Bbb Z}$-module, one has  $m_1g_1+m_2g_2+\cdots+m_kg_k \in G$
where $g_i\in G$ and $m_i\in {\Bbb Z}$. Note that when $m_i \geq
n$, then we can write $m_i = \ell n +r$ with $r<n$. Since $n$ is
the exponent of $G$, for any $g\in G$, we get, $m_ig = \ell ng +
rg = rg$. Also, if $m_i < 0,$ then $m_ig = |m_i|(-g).$ Therefore,
it is enough to vary the subset $A$ among the subsets of $]n[$
instead of  the set of all integers. Hence, among the relations of
the form $m_1g_1+m_2g_2+\cdots +m_kg_k$ with \mbox{$m_i\in\ ]n[$},
there may be   many  such linear  combinations equal to $0$ in
$G$. This motivates us to make the following definition.

\begin{deff}$\left.\right.$\vspace{.5pc}

\noindent Davenport's constant for $G$ with
respect to $A$ is denoted by $d_A(G)$ and is defined to be  the
least positive
integer $t$ such that given any sequence $S = (g_1, g_2, \dots,
g_t)$ in $G$,  we can always extract a subsequence $g_{i_1},
g_{i_2}, \dots, g_{i_\ell}$  with $1\leq  \ell \leq t$ such that\pagebreak
\begin{equation}
\sum_{j=1}^\ell a_{j}g_{i_j}  = 0 \mbox{ in } G,
\end{equation}
where $a_j \in A$ for all $j$.

When \hbox{$A = \{a\} \subset\ ]n[$} with $(a,n) = 1$, the
definition of $d_A(G)$ is nothing but the well-known Davenport's
constant which is denoted by $D(G)$. It is easy to prove that
$D({\Bbb Z}_n) = n$.

Olson, in \cite{ols1} and \cite{ols2},  proved that

\begin{enumerate}
\leftskip .1pc
\renewcommand{\labelenumi}{(\roman{enumi})}
\item when $G\sim {\Bbb Z}_m\oplus{\Bbb Z}_n$, where $1 < m|n$ integers,
we have $D(G) = m+n-1$;

\item when $G \sim {\Bbb Z}_{p^{e_1}}\oplus{\Bbb Z}_{p^{e_2}}\oplus
\cdots \oplus {\Bbb Z}_{p^{e_l}}$,  where $1\leq e_1 \leq e_2\leq
\cdots \leq e_l$ are integers, $D(G) = 1+\sum_{i=1}^l(p^{e_i}-1)$.
\end{enumerate}

It is seemingly a difficult problem to find the exact value of $D(G)$
for all $G$ other than the above mentioned groups.

If \hbox{$n\in A \subset\ ]n[$}, then, clearly, $d_A(G) = 1$.
Thus, we can always assume that \hbox{$A \subset\ ]n[$} and $1
\leq |A| < n$ and $n\not\in A$.
\end{deff}

\begin{probl}
Find the value of $d_A(G)$ for all non-empty
subsets $A$ of $]n[$ and for all finite abelian groups $G$ of exponent
$n$.

In this article, we  prove the following theorem.
\end{probl}

\begin{theor}[\!]
Let $G \sim {\Bbb Z}_{p^{e_1}}
\oplus{\Bbb Z}_{p^{e_2}}\oplus \cdots \oplus {\Bbb Z}_{p^{e_l}}$
where $1\leq e_1 \leq e_2\leq \cdots \leq e_l$ are integers.
Then{\rm ,} for any non-empty subset $A$ of $]p^{e_l}[$ such that
the elements of $A$ are incongruent modulo $p$ and non-zero modulo $p${\rm ,}
we have
\begin{equation*}
d_A(G) \leq \left\lceil\frac{1}{|A|}\left(1+\sum\limits_{i=1}^l(p^{e_i}-1)
\right)\right\rceil
\end{equation*}
where $\lceil x\rceil$ denotes the smallest positive integer greater
than or equal to $x$.
\end{theor}

\begin{coro}$\left.\right.$\vspace{.5pc}

\noindent Let $G \sim {\Bbb Z}_p^d,$ where $d\geq 1$ integer. Then
for any non-empty subset $A$ of $]p-1[${\rm ,} we have
\begin{equation*}
d_A(G) \leq \left\lceil\frac{1}{|A|}(d(p-1)+1)\right\rceil.
\end{equation*}
\end{coro}

\begin{coro}$\left.\right.$\vspace{.5pc}

\noindent Let $G \sim {\Bbb Z}_p^d,$ where $d\geq 1$ integer. If

\begin{enumerate}
\renewcommand{\labelenumi}{\rm (\alph{enumi})}
\leftskip .1pc
\item $A = \ ]p-1[\,\subset\ ]p[${\rm ,}  then  $d_A(G) = d+1${\rm ;}

\item $A_1 = \{a \in\ ]p-1[ : a \equiv x^2(\hbox{\rm mod}\ {p})\mbox{ for some }
x\in\ ]p[\}${\rm ,} then $d_{A_1}(G) =  2d+1${\rm ;}

\item $A_2 =\{a \in\ ]p-1[ : a \not\equiv x^2(\hbox{\rm mod}\ {p}) \mbox{ for all }
x\in\ ]p[\}${\rm ,} then $d_{A_2}(G) = 2d+1${\rm ;}

\item $A_3 =\{a \in\ ]p-1[ : a \mbox{ generates } {\Bbb F}_p^* \}${\rm ,}
then $d_{A_3}(G) = 2d+1${\rm ;}

\item $A_4 = \hbox{$\{a \in\ ]p-1[:$}\ a \in A_2\backslash A_3\},$ then
$d_{A_4}(G) = 2d+1$ holds for all primes $p\ne 2^{2^m}+1$ for some
$m\in{\Bbb N}.$

\item {$A_5 \subset\ ]p-1[$} such that $|A_5| = (p-1)/2$ and if $x\in
A_5${\rm ,} then $p-x\not\in A_5${\rm ,} then  $d_{A_5}({\Bbb Z}_p^d) = 2d+1.$

\item $A =\ ]r[\,\subset\ ]p[$ and $r\geq d${\rm ,} then
\begin{equation*}
\hskip -1.25pc d_A(G) = \left\lceil\frac{d(p-1)+1}{r}\right\rceil.
\end{equation*}
\end{enumerate}
\end{coro}

From Corollary 1.2(a), (b), (c), (f) and (g), it is clear that the
bound stated in Corollary~1.1 is tight for many subsets
\hbox{$A\subset\ ]p[$}, while for the subsets $A_3$ and $A_4$,
Corollary~1.1 provides a weaker bound, as $|A_3| = \phi(p-1)$ and
$|A_4| = (p-1)/2-\phi(p-1)$.\vspace{.5pc}

\noindent{\it Open problem.} For any finite abelian group $G$ of
exponent $n$, classify all the non-empty subsets $A$ of $]n[$ such
that
\begin{equation*}
d_A(G) \leq \left\lceil\frac{D(G)}{|A|}\right\rceil.
\end{equation*}
In  \S 3, we prove some elementary but sharp results for the
cyclic groups. In fact, one of the results says that when \hbox{$A
= \{a \in\ ]n[: (a,n) =1\}$}, then $d_A({\Bbb Z}_n) =
1+\Omega(n).$ It is easy to see that $1+\Omega(n) > d_A({\Bbb
Z}_n)/|A| = n/\phi(n)$ whenever $n = p^r$ for all integers $r\geq
r_0.$ Therefore, in the above open problem, the upper bound does
not hold for all non-empty  subsets $A$.

In \S 4, we explore the relation between $d_A(G)$ and
the associated constants for the existence of similar weighted sum
of length $|G|$ or the exponent of $G$. We prove some results and
state  a conjecture related to this relation.

\section{Proof of Theorem 1}

Let $p$ be any prime number. Throughout this section, we assume
that $G$ is a finite abelian $p$-group written {\it
multiplicatively}. Therefore, $G \sim {\Bbb Z}_{p^{e_1}}
\times{\Bbb Z}_{p^{e_2}}\times \cdots \times {\Bbb Z}_{p^{e_l}}$,
where $1\leq e_1 \leq e_2\leq \cdots \leq e_l$ are integers. Let
$n:= p^{e_l}$ be its exponent. Also, we denote the identity
element of $G$ by ${\bf e}$. We shall start with two lemmas which
are crucial for the proof of Theorem~1.

\setcounter{section}{1}
\begin{lem}\hskip -.3pc \hbox{\rm \cite{tro1}.}
Let $p$ be a prime number and $r\geq 1$ be an integer. Let $B =
\{0, b_1, b_2,\dots, b_r\}$ be a non-empty subset of
$]p-1[\,\cup\,\{0\}$ such that $b_i\not \equiv b_j(\hbox{\rm mod}\
{p})$ for all $i\ne j$. Then there exists  a polynomial
\begin{equation*}
f(x) = c_0+c_{b_1}x^{b_1}+ \cdots + c_{b_r}x^{b_r} \in {\Bbb F}_p[x]
\end{equation*}
such that $(1-x)^r$ divides $f(x)$ and $c_0 = f(0) \ne 0$.
\end{lem}

We recall that the group algebra ${\Bbb F}_p[G]$ is a ${\Bbb
F}_p$-vector space with $G$ as its basis. Hence, any element $\tau
\in {\Bbb F}_p[G]$ can be written as $\tau = \sum_{g\in G}a_g g$,
where $a_g\in {\Bbb F}_p$. Also, note that  $\tau = 0 \in {\Bbb
F}_p[G]$ if and only if $a_g = 0$ for all $g\in G$.

\begin{lem}\hskip -.3pc \hbox{\rm \cite{ols1}.}
 If $g_1, g_2, \dots,
g_s$ is a sequence in $G$ with $s\geq 1+\sum_{i=1}^l(p^{e_i}-1)${\rm ,}
then the element
\begin{equation*}
\prod\limits_{i=1}^s(1-g_i) = 0 \in {\Bbb F}_p[G].
\end{equation*}
\end{lem}

\begin{pot}
{\rm Let \hbox{$A = \{a_1,a_2,\dots,a_r\} \subset\ ]n[$} such that
$a_1 < a_2  < \dots < a_r$ and $a_i\nequiv a_j(\hbox{\rm mod}\
{p})$ for $i\ne j$. Set $B = A\cup \{0\}.$ Consider the group
algebra ${\Bbb F}_p[G]$.  Let $k\geq (1+
\sum_{i=1}^l(p^{e_i}-1))/r$ be any integer where $r = |A|.$ Let $S
= (g_1, g_2, \dots, g_k)$ be any given sequence in $G$ of length
$k$. To prove the theorem, it is enough to prove the existence of
a subsequence $g_{i_1}, \dots, g_{i_\ell}$ with $1\leq \ell\leq k$
such that
\begin{equation*}
\prod_{j=1}^\ell g_{i_j}^{a_{i_j}}  = {\bf e} \mbox{ in } G,
\end{equation*}
where  $a_{i_j} \in A$ for all $j$.

By Lemma 1.1, we have a polynomial $f(x) \in {\Bbb F}_p[x]$
associated with $B$ and $f(x) = (1-x)^rF(x)$ for some
polynomial $F(x) \in  {\Bbb F}_p[x]$.

Let $\sigma = \prod_{i=1}^k f(g_i)$. Clearly $\sigma$ is
the element in the group algebra ${\Bbb F}_p[G]$. In fact, as $G$
is abelian, we have
\begin{equation*}
\sigma =  \prod_{i=1}^k F(g_i)\prod_{i=1}^k(1-g_i)^r \in
{\Bbb F}_p[G].
\end{equation*}
Since $kr \geq 1+\sum_{i=1}^l(p^{e_i}-1)$, by
Lemma 1.2, we see that
\begin{equation*}
\prod_{i=1}^k(1-g_i)^r = 0 \mbox{ in } {\Bbb F}_p[G]
\end{equation*}
and hence, $\sigma = 0 \mbox{ in } {\Bbb F}_p[G]$.

Set $a_0=0$ and
\begin{equation*}
f(x) = \lambda_{a_0}+\lambda_{a_1}x^{a_1}+\cdots+\lambda_{a_r}
x^{a_r} \in {\Bbb F}_p[x].
\end{equation*}
Then, on the other hand, we can expand the product and see that
$\sigma$  is of  the following form:
\begin{equation}
0 = \sigma = \prod_{j=1}^k f(g_j) =
\prod_{j=1}^k\left(\sum_{i=0}^r\lambda_{a_i}g_j^{a_i}\right) \in
{\Bbb F}_p[G].
\end{equation}
The constant term of the above product is $\lambda_{a_0}^k {\bf e}$.
Since, by  Lemma 1.1, we know that  $\lambda_{a_0} \ne 0$, the
constant term in  (2) is non-zero. But since $\sigma = 0 \in
{\Bbb F}_p[G],$ we conclude that there is some other contribution to
${\bf e}$ in the above product. That is, we have
\begin{equation*}
\prod_{j=1}^\ell g_{i_j}^{a_{i_j}} = {\bf e} \mbox{ with }
a_{i_j} \in A.
\end{equation*}
Hence the theorem. $\hfill\Box$}
\end{pot}

\setcounter{theore}{0}
\begin{poc}
{\rm Since $G\sim {\Bbb Z}_p^d$, and \hbox{$A \subset\ ]p-1[$},
any two elements of $A$ are incongruent modulo $p$. Hence, by
Theorem 1, we get the result. $\hfill\Box$}
\end{poc}

\begin{poc}\hskip -.3pc \hbox{{\rm (a)}.} {\rm The upper bound
follows from  Corollary~1.1. Indeed, since \hbox{$A =\ ]p-1[,$} we
have $|A| = p-1.$  Therefore, by Corollary~1.1, we get
\begin{equation*}
d_A({\Bbb Z}_p^d) \leq \left\lceil\frac{d(p-1)+1}{p-1}\right\rceil
= d+1.
\end{equation*}
For the lower bound, consider the sequence $(e_1, e_2,\dots, e_d)$
in ${\Bbb Z}_p^d$ with $e_j = (0,0,\dots,0,1,0,\dots,0)$ where $1$
appears in the $j$th coordinate and $0$ elsewhere. Then clearly,
eq.~(1) is not satisfied for \hbox{$A =\ ]p-1[.$}

To  prove the assertions (b)--(e) at one stroke, we prove the following claim.

\begin{claim}
Let $S = (g_1, g_2, \dots, g_{2d+1})$ be any sequence in
${\Bbb Z}_p^d$ of length $2d+1$. Let $c\in {\Bbb Z}_p$ be a fixed non-zero element
such that
\begin{equation}
c \in \begin{cases}
A_1 & \mbox{ for proving  (b)}\\
A_2 & \mbox{ for proving  (c)}\\
A_3 & \mbox{ for proving  (d)}\\
A_4 & \mbox{ for proving  (e)}
\end{cases}.
\end{equation}
Then, for each $i = 1, 2, 3, 4$, the sequence $S$ satisfies  (1)
with coefficients $a_j$ in $A_i$ whenever $c\in A_i$.

First note that this claim clearly implies that $d_{A_i}
({\Bbb Z}_p^d) \leq 2d+1$ for all $i = 1, 2, 3, 4.$  Also note that
$A_4 = \emptyset$ whenever $p$ is of the form $2^{2^m}+1$. Hence
while proving (e), we need to assume that $p \ne 2^{2^m}+1.$

Since $g_j \in {\Bbb Z}_p^d$ for each $j = 1, 2, \dots, 2d+1$, we
put $g_j = (g_{1j},g_{2j},\dots, g_{dj})$
where $g_{lj}\in {\Bbb Z}_p$ for all $l = 1, 2, \dots, d$.
Let
\begin{equation*}
f_l(X_1,X_2,\dots,X_{2d+1}) = \sum_{j=1}^{2d+1}g_{lj}cX_j^2
\end{equation*}
for all $l = 1, 2, \dots\ \hbox{and}\ d$  be the system of homogeneous equations
over ${\Bbb Z}_p$. Since the total degree of $f_l$'s is equal to
$2d <$ the number  of variables involved, by the Chevalley--Warning
theorem, there exists a non-zero  solution in ${\Bbb Z}_p^{2d+1}$ to
the system. Let $(y_1,y_2,\dots,y_{2d+1}) \in {\Bbb Z}_p^{2d+1}$
be a non-zero solution and let
\begin{equation*}
I = \{j\in\{1,2,\dots, 2d+1\} : y_j\not\equiv 0(\hbox{\rm mod}\
{p})\},
\end{equation*}
which is non-empty. Therefore, we get
\begin{equation*}
0\equiv  \sum_{j\in I}g_{lj}cy_j^2(\hbox{\rm mod}\ {p}) \ \mbox{
for all } l = 1,2, \dots, d.
\end{equation*}
By putting $k_j = cy_j^2$ for all $j\in I$, we have
\begin{align*}
&\sum_{j\in I} g_{1j}k_j \equiv 0(\hbox{\rm mod}\ {p}), \sum_{j\in
I} g_{2j}k_j \equiv 0(\hbox{\rm mod}\ {p}), \dots, \sum_{j\in I}
g_{dj}k_j \equiv 0(\hbox{\rm mod}\ {p})
\end{align*}
and hence we arrive at
\begin{equation*}
\sum_{j\in I}g_jk_j = (0,0,\dots, 0) \mbox{ in } {\Bbb Z}_p^d.
\end{equation*}
Now, note that for each $i = 1, 2, 3, 4$, if $c\in A_i$, then
$k_j\in A_i$ for all $j \in I.$  Hence, the sequence $S$ satisfies
the claim.
\end{claim}\pagebreak

To prove the assertions (b)--(e), it is enough to prove the lower
bound.

Let $c$ be a fixed non-zero element in ${\Bbb Z}_p$ such that
\begin{equation*}
c = \Bigg\lbrace\begin{array}{cl}
\mbox{a quadratic non-residue modulo } p, & \mbox{ if } p\equiv 1(\hbox{\rm mod}\ {4})\\[.3pc]
1, & \mbox{ if } p\equiv 3(\hbox{\rm mod}\ {4})
\end{array}.
\end{equation*}
Consider the sequence  $S = (e_1, e_2,\dots, e_d, f_1, f_2, $
$\dots, f_d)$ in ${\Bbb Z}_p^d$  of length $2d$ where $e_j = (0,
0, \dots,0,1,0,\dots,0)$ and $f_j  = (0, 0,$ $\dots,0,
c,0,\dots,0)$ (in $e_j$ (similarly, in $f_j$), the element $1$
(similarly, $c$) appears in the $j$th coordinate  and $0$
elsewhere). If any  subsequence of  $S$ satisfies eq.~(1), then,
in the $j$th coordinate, we  have the following:
\begin{equation*}
a+cb \equiv 0(\hbox{\rm mod}\ {p}),\quad \mbox{ where } a, b\in
A_i.
\end{equation*}
Note that $a$ and $a^{-1}\in {\Bbb Z}_p^*$ are either both
quadratic residues or both  quadratic non-residues. Therefore,
$ab^{-1}$ is a quadratic residue modulo $p$, as $a, b \in A_i$. We
know that $p\equiv 1(\hbox{\rm mod}\ {4})$ if and only if  $-1$ is
a quadratic residue modulo $p$.  Therefore, we get $-ab^{-1}$ is a
quadratic residue  modulo $p$, whenever $p\equiv 1(\hbox{\rm mod}\
{4})$ and $-ab^{-1}$ is a quadratic non-residue, if $p\equiv
3(\hbox{\rm mod}\ {4})$. This contradicts the fact that $c\equiv
-ab^{-1}(\hbox{\rm mod}\ {p})$ is a quadratic non-residue modulo
$p$, in the case when $p\equiv 1(\hbox{\rm mod}\ {4})$, while if
$p\equiv 3(\hbox{\rm mod}\ {4})$, $c = 1$ is a quadratic residue.
Hence, for each $i = 1, 2, 3, 4$, the sequence $S$ does not
satisfy eq.~(1) with coefficients $a_j \in A_i$. Thus, we arrive
at the assertions (b)--(e).\vspace{.5pc}

\noindent (f) By assumption \hbox{$A_5 \subset\ ]p-1[$} and $|A_5|
= (p-1)/2$. Also, if $x\in A_5$, then $p-x\not\in A_5$. Therefore,
to get the lower bound,  consider the sequence $(e_1,e_1,
e_2,e_2,\dots, e_d,e_d)$ in ${\Bbb Z}_p^d$ of length $2d$ and
clearly, eq.~(1) is not satisfied for $A_5.$ Since $|A_5| =
(p-1)/2$, the upper bound, by Corollary~1.1, is $d_{A_5}(G) \leq
2d+1$ and thus the result.\vspace{.5pc}

\noindent (g) The upper bound follows from Corollary~1.1. To prove the lower
bound, consider the sequence
\begin{equation*}
S = ({\underbrace{e_1,\dots,e_1}_{m\mbox{ times }}}, \dots,
  {\underbrace{e_d,\dots,e_d}_{m\mbox{ times }}})
\end{equation*}
in ${\Bbb Z}_p^d$ length $dm$ where $m = \left\lceil\frac{p}{r}
\right\rceil -1 = \left[\frac{p}{r}\right].$ Clearly, any subsequence of $S$
does not satisfy  (1) and hence $d_A(G) \geq dm+1 = d\left[\frac{p}{r}
\right]+1$. Since $r\geq d$, we have
\begin{equation*}
d\left[\frac{p}{r} \right]+1 = \left\lceil\frac{d(p-1)+1}{r}
\right\rceil.
\end{equation*}
Thus, the corollary follows.$\hfill \Box$}
\end{poc}

\setcounter{section}{2}
\section{Elementary results for the cyclic group}

Let $n$ be any composite  positive integer.

\begin{theor}[\!]
For all \hbox{$a\in\ ]n[,$} we have
\begin{enumerate}
\leftskip .4pc
\renewcommand{\labelenumi}{\rm (\roman{enumi})}
\item $d_{\{a\}}({\Bbb Z}_n) = n/(a,n)$ where $(a,n)$ denotes
the  gcd of $n$ and $a${\rm ;}

\item $d_{\{a,n-a\}}({\Bbb Z}_n) = 1+ \lfloor\log_2n\rfloor${\rm ,}
whenever $(a,  n) = 1${\rm ;}

\item whenever \hbox{$A = \{a\hbox{\rm :}\ (a,n) = 1\} \subset\
]n[,$} we have{\rm ,} $d_{A}({\Bbb Z}_n) = 1+ \Omega(n)${\rm ,}
where $\Omega(n)$ denotes the number of prime power divisors $> 1$
of $n${\rm ;}

\item whenever \hbox{$A =\ ]n-1[\,\subset\ ]n[,$} we have{\rm ,}
$d_{A}({\Bbb Z}_n) =2$.
\end{enumerate}
\end{theor}

\begin{proof}$\left.\right.$

\begin{enumerate}
\leftskip .4pc
\renewcommand{\labelenumi}{(\roman{enumi})}
\item Whenever $(n,a) = 1$, the classical
Davenport constant $D({\Bbb Z}_n)$ and $d_{\{a\}}({\Bbb Z}_n)$ are
same and therefore the result follows easily. Hence, we assume that
$(n,a) = d >  1$. Let $S = (a_1, a_2, \dots, a_l)$ be
any sequence in ${\Bbb Z}_n$ of length $l = n/d.$

\hskip 1pc We have to find a subsequence of $S$ satisfying
equation (1). That is, we need to find  a subsequence, say,
$a_{i_{1}},a_{i_{2}}, \dots, a_{i_{r}}$ such that
\begin{equation*}
\hskip -1.25pc a\sum_{j=1}^ra_{i_{j}} \equiv 0(\hbox{\rm mod}\
{n}) \Longrightarrow \frac{a}{d} \sum_{j=1}^r a_{i_{j}} \equiv
0\mbox{ mod}\left(\frac{n}{d}\right).
\end{equation*}
Since $(a/d, n/d) = 1$, it is enough to find a subsequence of $S$
whose sum  is divisible by $n/d$. Since $l \geq n/d$, this is
possible by the classical  Davenport constant for the group ${\Bbb
Z}_{n/d}$. Thus, we have the  required upper bound. The lower
bound follows from the sequence $(a, a, \dots, a)$  where $a$
appears exactly $-1+n/(n,a)$ times.

\item Given that $A = \{a, n-a\}$ with $(a, n) = 1.$ To
prove the lower bound, let $s = \lfloor\log_2n\rfloor$.  Then
clearly, $2^{s}\leq n < 2^{s+1}.$  Consider the sequence
$S = (1, 2, 2^2,\dots,2^{s-1})$ in ${\Bbb Z}_n$ of length $s$.
Then any zero sum with coefficients in $A$ leads to an equation
of the type
\begin{equation*}
\hskip -1.25pc a(x-y)\equiv 0(\hbox{\rm mod}\ {n}),
\end{equation*}
where
\begin{equation*}
\hskip -1.25pc x = \sum_{i\in I}2^i\quad \mbox{ and }\quad y =
\sum_{j\in J}2^j
\end{equation*}
and $I\cup J$ is a partition of $\{1, 2, \dots, s\}$. Since $a$ is
coprime to $n$, we get that $x \equiv y(\hbox{\rm mod}\ {n})$ and
since both $x$ and $y$ are nonnegative integers smaller than $n$
we get $x = y$, which is impossible because of the uniqueness of
the binary expansion of a nonnegative integer.  Hence, the
sequence $S$  does not satisfy (1).

\hskip 1pc For the upper bound, let  $S = (a_1,a_2, \dots,a_s)$ be
any sequence in  ${\Bbb Z}_n$ of length $s =  1+ \lfloor\log_2n
\rfloor.$ Consider the set
\begin{equation*}
\hskip -1.25pc \sum(S) := \left\{\sum_{i\in I}aa_i\!\!: \ I
\subset \{1, 2,\dots, s\}\right\}.
\end{equation*}
Clearly, the set $\sum(S)$ contains $2^s$ elements. Since $n < 2^s$,
it  follows that there exist $I \ne J$ subsets of $\{1, 2,\dots,
s\}$ and satisfying
\begin{align*}
\hskip -1.25pc \sum_{i\in I}aa_i \!\equiv\!  \sum_{j\in
J}aa_j(\hbox{\rm mod}\ {n}) \Longrightarrow \sum_{i\in I}aa_i
\!+\! \sum_{i\in I}(n-a)a_i\equiv 0 (\hbox{\rm mod}\ {n}).
\end{align*}
Note that if $i\in I\cap J$, then, in the above congruence, we
have $aa_i + (n-a)a_i = na \equiv 0(\hbox{\rm mod}\ {n}).$ Hence,
we can assume that $I\cap J = \emptyset$ and this proves the upper
bound.

\item To see the lower bound, let $n = p_1p_2\cdots p_s$ where
$p_i$s are  prime divisors (not necessarily distinct) of $n$. The
sequence  $S = (1, p_1, p_1p_2,\dots, p_1p_2\dots
p_{s-1})$ in ${\Bbb Z}_n$  of  length $\Omega(n)$ does not satisfy
eq.~(1).\pagebreak

For the upper bound\footnote{One can also give a direct and
elementary proof.}, we use a result of Luca \cite{luca} stated in
the next section. Let $k = 1+\Omega(n)$, $g_1, g_2, \dots, g_k$ be
any elements in ${\Bbb Z}_n$ and consider the sequence $S = (g_1,
\dots, g_k, \underbrace{0, \dots, 0}_{n-1\mbox{ times}})$ in
${\Bbb Z}_n$ of length $n+ \Omega(n)$. By Luca's result, there is
a linear combination with coefficients in  $A$ of precisely $n$
elements from  $S$ which is $0$. Of those $n$ elements, at most
$n-1$ elements are  $0$'s. Hence, we have zero sum of $g_i$'s with
coefficients in $A$ which proves the upper bound.

\item The lower bound follows by taking $k = 1$ and $g_1 = 1$.
For the upper bound, let $g_1, g_2$ be any two elements in $]n[$.
If one of them is $n$, say $g_1 = n$, then $1\cdot g_1 = 0$. If
both are $< n$, then $(n-g_2)g_1+g_1\cdot g_2 \equiv 0(\hbox{\rm
mod}\ {n})$ so we can take $a_1 = n-g_2$ and $a_2 = g_1$, both in
$]n-1[$, which proves the result. Hence, the corollary follows.
$\hfill\Box$
\end{enumerate}
\end{proof}

\vspace{-2pc}
\section{Relation between $\pmb{d_A(G)}$ and zero-sums of length $\pmb{|G|}$}

Let $G$ be a finite abelian group of exponent $n$ and $A$ be any
non-empty subset of $]n[$.

By $ZS_A(G)$ (similarly, $s_A(G)$), we
denote the least positive  integer $t$ such that for any given
sequence $S = (g_1, g_2, \dots, g_t)$ in $G$ of length $t$ has a
non-empty subsequence, say,  $g_{i_1},g_{i_2}, \dots,g_{i_m}$
satisfying
\begin{equation}
\sum_{j=1}^m a_jg_{i_j} = 0 \mbox{ in } G,
\end{equation}
where $a_j \in A$ for all $j$ and $m = |G|$ (similarly, $m = n$).

The following results are known for various subsets $A$ and some
groups:
\begin{enumerate}
\renewcommand{\labelenumi}{\arabic{enumi}.}
\item When $A = \{a\}$ and $(a,n) =1$, we have
\begin{equation*}
\hskip -1.25pc ZS_A(G) = D(G) +|G| -1
\end{equation*}
by \cite{gao1},  which generalizes the famous theorem of
Erd\"{o}s, Ginzburg and Ziv \cite{egz}.

\item When \hbox{$A = \{a, n-a\}\subset\ ]n[$} and $(a,n) = 1$, we have
\begin{equation*}
\hskip -1.25pc s_A({\Bbb Z}_n) = ZS_A({\Bbb Z}_n) =
n+\lceil\log_2n\rceil
\end{equation*}
by \cite{adhi1}.

\item When \hbox{$A = \{a\in\ ]n[ : (n,a) =1\}$}, we have
\begin{equation*}
\hskip -1.25pc s_A({\Bbb Z}_n)= ZS_A({\Bbb Z}_n) = n+\Omega(n)
\end{equation*}
by \cite{luca}.

\item When \hbox{$A =\ ]r[\,\subset\ ]p[$}, we have
\begin{equation*}
\hskip -1.25pc s_A({\Bbb Z}_p) = ZS_A({\Bbb Z}_p) = p+[p/r]
\end{equation*}
by \cite{adhi2}.
\end{enumerate}

Note that, by the definition, $s_A({\Bbb Z}_n) = ZS_A({\Bbb
Z}_n)$. Also, we have $ZS_A(G) > d_A(G) +|G| -2$. Indeed,  by the
definition of $d_A(G)$, we have a sequence $S$ in $G$ of length
$d_A(G)-1$ and does not satisfy (1). Now consider the sequence
\begin{equation*}
T = (S, \underbrace{0,0,\dots,0}_{|G|-1 \ {\rm times}})
\end{equation*}\pagebreak

\noindent in $G$ of length $d_A(G)+|G|-2$. Clearly, by the
construction, $T$ does not satisfy (4) with $m = |G|.$ Therefore,
$ZS_A(G)\geq d_A(G) +|G| -1$. We feel that the lower bound seems
to be tight. More precisely, we have Theorems 3 and 4.

\begin{conjec}
{\rm For any finite abelian group $G$
with  exponent $n$ and for any non-empty subset $A$ of $]n[$, we
have
\begin{equation*}
ZS_A(G) = |G|-1+d_A(G).
\end{equation*}
Using Theorem 2 and the  results 1--4 stated above, we see that
Conjecture 1 holds for those $A$'s and $G$'s. In support of
Conjecture~1, we prove Theorems~3 and 4.}
\end{conjec}

\setcounter{section}{3}
\setcounter{defin}{0}
\begin{lem}
If $s_A({\Bbb Z}_n^d) = d_A({\Bbb Z}_n^d) +n-1$ for some $A
\subset\ ]n[,$ then we have
\begin{equation*}
ZS_A({\Bbb Z}_n^d) = d_A({\Bbb Z}_n^d) +n^d-1.
\end{equation*}
\end{lem}

\begin{proof}
Consider a sequence
$S = (g_1,g_2, \dots, g_\ell)$ in ${\Bbb Z}_n^d$ of length $\ell =
d_A({\Bbb Z}_n^d) +n^d-1$. To prove the lemma, we have to prove that
$S$ has a subsequence satisfying (4) with $m = n^d$. By hypothesis,
we know that
$s_A({\Bbb Z}_n^d) = d_A({\Bbb Z}_n^d) +n-1$. Since $\ell \geq
s_A({\Bbb Z}_n^d)$, clearly, there exists a subsequence $S_1$ of $S$
of length $n$ satisfying (4) with $m = n$. Since
\begin{equation*}
|S| = \ell =  d_A({\Bbb Z}_n^d) +n^d-1 =
d_A({\Bbb Z}_n^d) +n(n^{d-1}-1)+n-1,
\end{equation*}
we can extract disjoint subsequences $S_1, S_2, \dots,
S_k$ of $S$ with $|S_i| = n$ for all $i = 1,2,\dots,k$ where $k =
n^{d-1}$ satisfies (4). Note that the total length of the
subsequence $S_i$ is $n^{d-1}n = n^d$. Thus, we get
\begin{equation*}
\sum_{i=1}^k\sum_{j=1}^na_{ij}g_{ij} \equiv 0(\hbox{\rm mod}\
{n}),
\end{equation*}
where $a_{ij}\in A$ and $g_{ij}\in S_i$ for all $i = 1,2,\dots,k$
and for all $j=1,2,\dots,n$. Hence, we arrive at a subsequence $L$
of $S$ of length $n^d$ satisfying (4). $\hfill\Box$
\end{proof}

\begin{theor}[\!]
Let $d\geq 1$ be any
integer. Let $p$ be a prime number such that $p\geq 2d+1$. Then
\begin{equation*}
s_{A_i}({\Bbb Z}_p^d) = d_{A_i}({\Bbb Z}_p^d)+ p-1 = 2d+p,
\end{equation*}
where $A_i$\hbox{\rm '}s {\rm (}for $i= 1, 2, 3, 4${\rm )} are as
defined in Corollary {\rm 1.2.}
\end{theor}

\begin{proof}
Choose \hbox{$c\in\ ]p-1[$} as in (3). Let $S =
(g_1,g_2,\dots,g_{2d+p})$ be a  given sequence in ${\Bbb Z}_p^d$
of length $2d+p$. To conclude the proof of this theorem,  we have
to prove that for each $i = 1,2,3,4$ the sequence $S$ has a
subsequence which satisfies (4) for $A_i$  with $m = p$.

Since $g_i\in {\Bbb Z}_p^d$, we have
\begin{align*}
g_i = (g_{i1}, g_{i2}, \dots, g_{id}) \mbox{ where } g_{ij}\in
{\Bbb Z}_p \mbox{ for all } j = 1, 2, \dots, d.
\end{align*}
Now consider the homogeneous equations over the finite field
${\Bbb F}_p$ as follows:
\begin{equation*}
f_j(X_1,X_2, \dots, X_{2d+p}) = \sum_{i=1}^{2d+p}cg_{ij}X_i^2\quad
\mbox{ for all } j = 1,2,\dots, d
\end{equation*}
and
\begin{equation*}
g(X_1,X_2,\dots, X_{2d+p}) = \sum_{i=1}^{2d+p} X_i^{p-1}.
\end{equation*}
Note that sum of the degrees of $f_j$ and $g$ is $2d+p-1$ which
is strictly less than the number of variables $X_i$'s. Hence, by the
Chevalley--Warning theorem, we have a non-trivial simultaneous
solution over ${\Bbb F}_p$.  Let $(y_1,y_2,\dots,y_{2d+p}) \in
{\Bbb F}_p^{2d+p}$ be a non-trivial solution and let
\begin{equation*}
I = \{i\in\{1,2,\dots, 2d+p\}\!: y_i\not\equiv 0(\hbox{\rm mod}\
{p})\} \neq \emptyset.
\end{equation*}
Then, we get
\begin{equation*}
0\equiv  \sum_{i\in I}cg_{ij}y_i^2 (\hbox{\rm mod}\ {p})
\end{equation*}
and
\begin{equation*}
0\equiv  \sum_{i\in I}y_i^{p-1} \equiv \sum_{i\in I} 1= |I|
(\hbox{\rm mod}\ {p}),
\end{equation*}
as $y_i\not\equiv 0(\hbox{\rm mod}\ {p})$ whenever $i\in I$ and by
Fermat Little theorem, the above follows. Note that $|I| \ne 0$
and $|I| \leq 2d+p < 2p$. Thus, using this fact and the second
congruence above, we get $|I| = p$. From the first congruence, we
get
\begin{equation*}
\sum_{i\in I}g_icy_i^2 \equiv 0(\hbox{\rm mod}\ {p}),
\end{equation*}
where $(g_i)_{i\in I}$ is a subsequence of the given sequence with
$|I| = p$. Also, note that $cy_i^2 \in A_j$ for all $i \in I$,
whenever $c\in A_j$. Hence, the theorem follows. $\hfill\Box$
\end{proof}

\setcounter{theore}{0}
\begin{coro}$\left.\right.$\vspace{.5pc}

\noindent Let $d\geq 1$ be an
integer and $p\geq 2d+1$ be any prime number. Then
\begin{equation*}
ZS_{A_i}({\Bbb Z}_p^d) = d_{A_i}({\Bbb Z}_p^d) + p^d-1= p^d+2d,
\end{equation*}
for all $A_i$\hbox{\rm '}s where $i = 1,2,3,4$ as defined in
Corollary {\rm 1.2(b), (c), (d)} and {\rm (e)}.
\end{coro}

\begin{proof}
Proof follows from Lemma 3.1 and Theorem 3.
$\hfill\Box$
\end{proof}

\begin{coro}$\left.\right.$\vspace{.5pc}

\noindent Let $d\geq 1$ be any integer. Let $p\equiv 3(\hbox{\rm
mod}\ {4})$ be a prime such that $p\geq 2d+1.$ Then
\begin{equation*}
ZS_{A_5}(\Bbb Z_p^d) = d_{A_5}({\Bbb Z}_p^d) +
p^d-1 = 2d+p^d,
\end{equation*}
where $A_5$ is defined in Corollary {\rm 1.2(f)}.
\end{coro}

\begin{proof}
Since $p\equiv 3 (\hbox{\rm mod}\ {4})$, we know that $-1 = p-1$
is a quadratic non-residue modulo $p$. Therefore, $x$ is a
quadratic residue modulo $p$, $p-x$ is a quadratic non-residue
modulo $p$ and vice versa. Thus, we can re-write
\begin{equation*}
A_5 = \{a\in\ ]p-1[ \ : \ a\equiv x^2(\hbox{\rm mod}\ {p})\mbox{
for some } x\in\ ]p-1[\},
\end{equation*}
or its compliment in $]p-1[.$ Hence, the result follows by
Corollary 3.1. $\hfill\Box$
\end{proof}

\begin{theor}[\!]
Let $S =  (g_1,g_2,\dots,g_k)$ be a  sequence in  $G$ of length $k
= |G| +$ $d_A(G)-1$. If $0 \in G$ appears in $S$  at least
$d_A(G)-1$ times{\rm ,} then $S$ satisfies eq.~{\rm (4)}
with\break $m = |G|$.
\end{theor}

\begin{proof}
By rearranging the terms of $S$, we may assume that
\begin{equation*}
S = ({\underbrace{0,0,\dots,0}_{h \ {\rm times}}}, g_1,g_2,\dots,
g_{k-h}),
\end{equation*}
where $g_i\ne 0 \in G$. If $h\geq |G|$, then eq.~(4) is
trivially satisfied. Assume that $h\leq |G|-1$. Therefore, we have
\begin{equation*}
k-h = |G| +d_A(G)-1-h \geq d_A(G).
\end{equation*}
Therefore, by the definition of $d_A(G)$, we can find $W_1 =
(g_{i_{1}}, g_{i_{2}}, \dots, g_{i_{r}})$, a subsequence of $S$
with $g_{i_{j}}\neq 0$ and $a_{i_{j}} \in A$ for $j = 1, 2, \dots,
r$ satisfying eq.~(1) with $r \geq 1$. Choose $W_1$ to be the
maximal subsequence having the above property. If $k-h-|W_1|\geq
d_A(G)$, then again we  choose another subsequence $W_2$ with
non-zero elements $g_i$s  satisfying eq.~(1). Hence, $W_1W_2$
together satisfy eq.~(1) which contradicts the maximality of
$W_1$. This forces that $k-h-|W_1| \leq d_A(G)-1$.  Hence, we
arrive at $|W_1| \geq |G|-h$.\break Therefore,
\begin{equation*}
|G|-h \leq |W_1|\leq k-h \leq |G|.
\end{equation*}
Since we have at least $h$ number of 0's outside $W_1$,
eq.~(4) is satisfied with $m = |G|$. Thus, the result is
proved. $\hfill\Box$
\end{proof}

\section*{Acknowledgment}

The author is grateful to the
referee for  his/her extensive comments to improve the presentation
of the paper.

\noindent{\it Note added.} It seems that much before H~Davenport
introduced the problem of \hbox{Davenport}'s constant in 1966 (see
for instance, H Davenport, Proceedings of the Midwestern
\hbox{Conference} on Group Theory and Number Theory, {Ohio State
University,} April, 1966), the problem was, historically, first
studied and introduced by K~Rogers in 1962 (K~Rogers,  A
combinatorial problem in abelian groups, {\it Proc. Cambridge
Philos.  Soc.} {\bf 59} (1963) \hbox{559--562}). Somehow this
reference was overlooked and never quoted by the later\break
authors.

\end{document}